\documentclass[%
 jmp,
 amsmath,amssymb,
notitlepage]{revtex4-1}

\usepackage{graphicx}
\usepackage{dcolumn}
\usepackage{bm}

\usepackage[utf8]{inputenc}
\usepackage[T1]{fontenc}
\usepackage{mathptmx}
\newcommand{\beq}{\begin{equation}}
\newcommand{\eeq}{\end{equation}}
\newcommand{\bea}{\begin{eqnarray}}
\newcommand{\eea}{\end{eqnarray}}
\newcommand{\bwd}{\begin{widetext}}
\newcommand{\ewd}{\end{widetext}}

\begin{document}

\preprint{AIP/123-QED}

\title[Sample title]{A few properties of the ratio of Davenport-Heilbronn Functions}

\author{Tao Liu}
\affiliation{Coherent Electron Quantum Optics Research Center, State Key Laboratory of Environment-friendly Energy Materials, Southwest University of Science and Technology, 59 Qinglong Road, Mianyang, Sichuan 621010, China\\
School of Science, Southwest University of Science and Technology, Mianyang, Sichuan 621010, China
}%

\author{Juhao Wu}
\affiliation{%
Stanford University, Stanford, California 94309, USA
}%

\date{\today}

\begin{abstract}
Starting from the Davenport-Heilbronn function equation: $f(s) = X(s) f(1-s)$, we discover the four properties of the meromorphic function $X(s)$ defined as the ratio of the Davenport-Heilbronn functions: $\frac{f(s)}{f(1-s)} = X(s)$, and three corresponding lemmas. For the first time, we propose to study the distribution of the non-trivial zeros of the Davenport-Heilbronn function by exploring the monotonicity of the similarity ratio $\left| \frac{f(s)}{f(1-s)} \right|$. We point out that for the $f(s)$ which satisfies the Davenport-Heilbronn function equation, the existence of non-trivial zeros outside of the critical line $\{ s_n | \sigma \neq 1/2 \}$ presents two puzzles: 1) $f(s_n) \neq f(1 - s_n)$; 2) the existence of non-trivial zeros $\{ s_n | \sigma \neq 1/2 \}$ is in contradiction of the monotonicity of the similar ratio $\left| \frac{f(s)}{f(1-s)}\right|$.
\end{abstract}

\maketitle

\section{Equation for Davenport-Heilbronn function}
The Davenport-Heilbronn function satisfies the function equation:
    \beq\label{zeta}
    f(s) = \left(\frac5{\pi}\right)^{\frac12 -s}\frac{\Gamma\left(1-\frac s2\right)}{\Gamma\left(\frac{1+s}2\right)}f(1-s),
    \eeq
where $\Gamma(1-s)$ is the analytical continuation of the factorial and $s = \sigma + i t$, with $\sigma \in R$ and $t \in R$ both being real number. Notice that $s = - 2 n - 1$ ($n = 1, 2,3, \cdots, \infty$) are the trivial zeroes of $f(s)$. Since $f(s)$ is a Holomorphic function, $f(s)$ is analytic everywhere on the complex $s$-plane.

\section{The Meromorphic Function $X(s)$}\label{Riemann}

\subsection{Introduction of the meromorphic function $X(s)$}
Based on Davenport-Heilbronn function equation (\ref{zeta}), we can introduce a meromorphic function:
    \beq\label{xs}
    X(s) = \left(\frac5{\pi}\right)^{\frac12 -s}\frac{\Gamma\left(1-\frac s2\right)}{\Gamma\left(\frac{1+s}2\right)},
    \eeq
so that Eq. (\ref{zeta}) is rewritten as:
    \beq\label{zetawz}
    f(s) = X(s) f(1-s).
    \eeq
The distribution of the nontrivial zeros of the function $f(s)$ satisfying Eq. (\ref{zeta}) is closely related to the properties of the meromorphic function $X(s)$. In the following, let us discuss the properties of $X(s)$.

\subsection{The properties of the meromorphic function $X(s)$}

\subsubsection{The reflection symmetry of $X(s)$}\label{property1} In the $s$-complex space, for arbitrary $\varepsilon \in R$, setting $s_0 = 1/2 + it$, then the pair: $s_{\pm} = s_0 \pm \varepsilon$ is a mirror symmetric pair with respect to $s_0$.
The complex conjugates of $s_+$ and $s_-$ are noted as $s^{\ast}_+$ and $s^{\ast}_-$. Under the $X(s)$ map, $s_{\pm}$ and their mirror reflected complex conjugate $s^{\ast}_{\mp}$ are reciprocal pair.

    \beq\label{ww1}
    X\left(s_{\pm}\right)X\left(s^{\ast}_{\mp}\right) \equiv 1.
    \eeq

\noindent {\it Proof}

On the complex plane of $s$, for arbitrary real number $\varepsilon \in R$, because:

    $$
    X(s_{+})
    = \left(\frac5{\pi}\right)^{-\varepsilon - i t} \frac{\Gamma\left(\frac34 - \frac{\varepsilon}2 - \frac{it}2\right)}{\Gamma\left(\frac 34 + \frac{\varepsilon}2 + \frac{it}2\right)},
    $$
and
    $$
    X(s^{\ast}_{-})
    =
    \left(\frac5{\pi}\right)^{\varepsilon + i t}\frac{\Gamma\left(\frac34 + \frac{\varepsilon}2 + \frac{it}2\right)}{\Gamma\left(\frac 34 - \frac{\varepsilon}2 - \frac{it}2\right)} = \frac 1{X(s_+)},
    $$
which is to say that in the $X$-space, $X(s^{\ast}_-)$ must be equal to $1/X(s_{+})$, which is the reciprocal of $X(s_+)$; we have:
    $$
    X(s_+)X(s^{\ast}_-) = X(s_+) \frac1{X(s_+)} = 1.
    $$
Similarly, we can prove that in the $X$-space, $X(s^{\ast}_+)$ must be equal to $1/X(s_{-})$, which is the reciprocal of $X(s_-)$, {\it i.e.}, $X(s^{\ast}_+)X(s_-) = 1$. Therefore, property \ref{property1} is proven.

\subsubsection{$X(s)$ has trivial zeros at $s = \{- 2n -1 | n = 0, 1,2, \cdots, \infty\}$ and pole $s = \{2 n + 2 | n = 0, 1, 2, \cdots, \infty \}$. $X(s)$ at these trivial zeros and at these poles are reciprocal pairs: $X(-2n-1)X(2 + 2n) \equiv 1$.}\label{reflective}\hfill

\noindent {\it Proof}

First of all, it is obvious that:
$$
X(-2n-1) = \left( \frac5{\pi}\right)^{\frac32+2n}\frac{\Gamma\left(\frac32+n\right)} {\Gamma(-n)} = 0
$$
for ($n=0,1,2,\cdots$),
and
$$
X(2n+2) = \left( \frac5{\pi}\right)^{-\frac32-2n}\frac{\Gamma(-n)}{\Gamma\left(\frac32+n\right)} = \infty$$
for ($n=0,1,2,\cdots$).

Apparently, $X(-2n-1)$ and $X(2n + 2)$ satisfy reflection symmetry relation:
    \beq
    X(-2n-1)X(2n + 2) = 1.
    \eeq
This then serves as the proof for property \ref{reflective}.

\subsubsection{Monotonicity of the absolute value $|X(s)|$ of $X(s)$. \newline \indent Excluding the zeros and the poles, the absolute value $|X(s)|$ of the map $X(s)$ is a monotonic function of $t$ except when $\sigma = 1/2$: \newline \indent 1. in the range $0 < t < +\infty$, when $\sigma > 1/2$, $|X(s)|$ monotonically decreases with the increase of $t$; when $\sigma < 1/2$, $|X(s)|$ monotonically increases with the increase of $t$. \newline \indent 2. in the range $-\infty < t < 0$, when $\sigma > 1/2$, $|X(s)|$ monotonically increases with the increase of $t$; when $\sigma < 1/2$, $|X(s)|$ monotonically decreases with the increase of $t$.}\label{mono}\hfill

\noindent {\it Proof}
Based on the series representation of the Digamma Function $\Psi(s)$ \cite{Liu20}, we have:

When $t > 0$,
    \beq\label{tg0}
    \frac{\partial}{\partial t}\frac{|f(s)|}{|f(1-s)|} = \frac{\partial}{\partial t} |X(s)| = \sum^{\infty}_{n=1}\frac{8t\left(\frac12-\sigma\right)\left(n-\frac14\right)|X(s)|}{| it + 2n + \sigma -1 |^2 | it + 2 n - \sigma |^2}
    \left\{
    \begin{array}{l}
    >0, \left( \sigma < \frac12 \right) \\
    =0, \left( \sigma = \frac12 \right) \\
    <0, \left( \sigma > \frac12 \right)
    \end{array};
    \right.
    \eeq
when $t < 0$,
    \beq
    \frac{\partial}{\partial t}\frac{|f(s)|}{|f(1-s)|} = \frac{\partial}{\partial t} |X(s)| = \sum^{\infty}_{n=1}\frac{8t\left(\frac12-\sigma\right)\left(n-\frac14\right)|X(s)|}{| it + 2n + \sigma -1 |^2 | it + 2 n - \sigma |^2}
    \left\{
    \begin{array}{l}
    <0, \left( \sigma < \frac12 \right) \\
    =0, \left( \sigma = \frac12 \right) \\
    >0, \left( \sigma > \frac12 \right)
    \end{array}.
    \right.
    \eeq

Therefore,
\newline \indent in the range $0 < t < +\infty$, when $\sigma > 1/2$, $|X(s)|$ monotonically decreases with the increase of $t$; when $\sigma < 1/2$, $|X(s)|$ monotonically increases with the increase of $t$.
\newline \indent in the range $-\infty < t < 0$, when $\sigma > 1/2$, $|X(s)|$ monotonically increases with the increase of $t$; when $\sigma < 1/2$, $|X(s)|$ monotonically decreases with the increase of $t$.

\noindent So, property \ref{mono} is proven.

\subsubsection{On the complex $s$-plane, in the range: $0 \leq \sigma \leq 1$, but $\sigma  \neq 1/2$, the $t$ satisfying $|X(s)|=1$ is bounded: $|t| < \kappa$.}\label{bounded}\hfill

\noindent {\it Proof}

1) Because of property \ref{mono}: {\it i.e.}, the monotonicity of $|X(s)|$, in the range $\sigma < 1/2$, for arbitrary $\sigma$, $|X(s)|$ always monotonically increases with the increase of $|t|$. Therefore, the $t$ satisfying $|X(s)| = 1$ must be bounded.

2) According to property \ref{property1}: due to the reflection symmetry of $X(s)$, it must be true that $|X(s_+)||X(s_-)| = 1$, {\it i.e.}, when $s_{\pm} = \frac 12 \pm \varepsilon + it$, $|X(s_+)|$ and $|X(s_-)|$ have reflection symmetry. So, in the range of $\sigma > 1/2$, for an arbitrary $\sigma$, $|X(s)|$ monotonically decreases when $|t|$ increases. Therefore, the $t$ satisfying $|X(s)| = 1$ must be bounded.

\begin{figure}[!htp]
 \centering
 \includegraphics[width=6cm,angle=0]{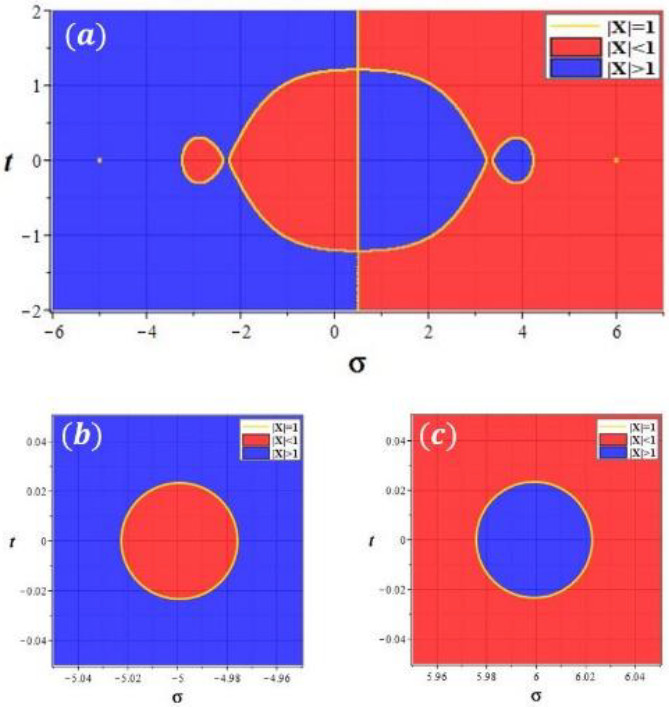}
 \caption{\label{ImpFunc}
 The implicit function curve $|X(s)| = 1$ on the $s$-plane for $-6 \leq \sigma \leq 7$: the yellow solid circle represents $|X(\sigma)| = 1$, the blue region represents $|X(s)| > 1$, and the red region represents $|X(s)| < 1$. The (b) and (c) are the zoom in details of (a) around the zero ($s = -5$) and the pole ($s = 6$) of $|X(s)|$.
 }
\end{figure}

3) The implicit function curve on the $s$-plane for $|X(s)| = 1$ is shown as in Fig. \ref{ImpFunc}.

In Fig. \ref{ImpFunc}, we have a clear illustration of the distribution of $|X(s)| = 1$ on the $s$-plane: on both sides of $\sigma = 1/2$, it encloses all the zeros and poles of $X(s)$, and they are symmetrically distributed. In the range $0 \leq \sigma \leq 1$, the $t$ satisfying $|X(s)| = 1$ is bounded $|t| < \kappa = 1.21164$.

\noindent So, property \ref{bounded} is proven.

\section{The distribution of the nontrivial zeros of the function $f(s)$}\label{Proof}
Let us take the absolute value of the Davenport-Heilbronn function equation in (\ref{zetawz}) to have:
    \beq
    |f(s)| = |X(s)| |f(1-s)|.
    \eeq
In the following, we will first prove a few lemmas of the function $f(s)$ based on the properties of the meromorphic function $X(s)$.

{\bf Lemma 1: functions $f(s)$ and $f(1-s)$ do not have nontrivial zeroes in the range: $0 < |X(s)| \neq 1$.}

\noindent {\it Proof}

{\bf 1) If $|X(s)| \neq 1$, then it must be true that $|f(s)| \neq |f(1 - s)|$.}

Assuming that $|f(s)| = |f(1-s)|$, based on the analytical properties of $f(s)$, it can be proven that $|X(s)| = 1$ (please refer to Appendix \ref{equiv}). Then this is in conflict of $|X(s)| \neq 1$. Therefore, the assumption could not hold and it has to be true that $|f(s)| \neq |f(1-s)|$.

{\bf 2) When $|X(s)| > 0$, for the function satisfying $|f(s)| \neq |f(1-s)|$, it is impossible to have $|f(s)| = 0$.}

Assuming that $|f(s)| = 0$, {\it i.e.}, $s$ is a nontrivial zero of $f(s)$, we would have $|X(s)| \neq 0$. Then due to $|f(s)| = |X(s)| |f(1-s)|$, it must be true that $|f(1-s)| = 0$, so that $|f(s)| = |f(1-s)|$. Then this is in conflict with the statement that $|f(s)| \neq |f(1-s)|$. Therefore, the assumption $f(s) = 0$ can not hold.

{\bf 3) When $|W(s)| > 0$, for the function satisfying $|f(s)| \neq |f(1-s)|$, $|f(1 - s)| = 0$ can only be valid on the trivial zeroes.}

Assuming that $|f(1 - s)| = 0$, then due to $|f(s)| = |X(s)| |f(1-s)|$, one would have $|f(s)| = 0$ with the only exception of $|X(s)| = \infty$, {\it i.e.}, $s$ is a pole of $|X(s)|$. However, $|f(s)| = 0$ is in conflict with $|f(s)| \neq |f(1-s)|$. Now, the poles of $|X(s)|$ are $s = 2 n + 2$, where $|f(1-s)| = |f(-2n -1)| = 0$, {\it i.e.}, these are just the trivial zeroes.

Based on the above 1), 2), and 3), we know that the function does not have nontrivial zeroes for $0 < |X(s)| \neq 1$. So Lemma $1$ is proven.

{\bf Lemma 2: $|X(s)|=1$ is the necessary condition for the nontrivial zeroes of the function $f(s)$.}\label{w1necess}

\noindent {\it Proof}

Based on Lemma $1$: In the range: $0 < |X(s)| \neq 1$, there is no nontrivial zeroes of the function $f(s)$.

Therefore, the nontrivial zeroes of the function $f(s)$ can only be on $|X(s)|=1$. Indeed, when $|X(s)| = 1$, we must have:
    \beq
    |f(s)|=|f(1-s)|.
    \eeq
If $|f(s)|=0$, then $|f(1-s)|=0$. However, $|f(s)|=|f(1-s)|$ does not guarantee $f(s)$ to be zero. Therefore, $|X(s)| = 1$ is only the necessary condition of the nontrivial zeroes of the function $f(s)$.

So Lemma $2$ is proven.

{\bf Corollary 1: for $s=1/2+it$ being the nontrivial zeroes of $f(s)$, the necessary condition for $X(s)$ is that $|X(s)|=1$.}

Inserting $s = \frac12 + i t$ and $s^{\ast} = \frac12 - i t$ into Eq. (\ref{xs}), we have:
    \beq
    X(s){\overline X(s)} = \left(\frac5{\pi}\right)^{1-s-s^{\ast}} \left.\frac{\Gamma\left(1-\frac s2\right)\Gamma\left(1-\frac{s^{\ast}}2\right)}{\Gamma\left(\frac{1+ s}2\right) \Gamma\left(\frac {1+s^{\ast}}2\right)}\right|_{s=\frac12+it,s^{\ast}=\frac12-it} = 1.
    \eeq
So we have: $|X(s)|=1$. The geometric illustration is that $X(s)$ maps the straight line $s=1/2+it$ in $s$-space to the unit circle in the $X$-space.

{\bf Lemma 3: on the complex $s$-plane, except for $s = \frac12 + it$, the similarity ratio $\frac{f(s)}{|f(1-s)|}$ of the function $|f(s)|$ does not have the form of $\frac 00$ for finite $t$.}

\noindent {\it Proof}

{\bf 1) for $\sigma \neq \frac12$, when and only when $t \rightarrow \pm \infty$, $\frac{f(s)}{|f(1-s)|}$ is in the form of $\frac 00$.}

According to Eq. (\ref{xs}) and Eq. (\ref{zetawz}), we have:
    \beq\label{absxs}
    |X(s)| = \frac{f(s)}{|f(1-s)|} = \left( \frac5{\pi} \right)^{\frac12 - \sigma} \frac{\left| \Gamma\left(1 - \frac s2\right) \right|}{\left| \Gamma\left( \frac{1+s}2 \right) \right|}.
    \eeq
Therefore:
    \beq
    \frac d{dt}\left| \Gamma\left( 1 - \frac s2 \right) \right| = -t \sum^{\infty}_{n = 1} \frac{\left| \Gamma\left(1 - \frac s2\right) \right|}{|\sigma + it - 2n|^2}
    \left\{
    \begin{array}{cc}
    < 0 & (t > 0) \\
    > 0 & (t < 0)
    \end{array}
    \right.,
    \eeq
    \beq
    \frac d{dt}\left| \Gamma\left( \frac{1 + s}2 \right) \right| = -t \sum^{\infty}_{n = 1} \frac{\left| \Gamma\left(\frac{1 +  s}2\right) \right|}{|\sigma + it + 2n - 1|^2}
    \left\{
    \begin{array}{cc}
    < 0 & (t > 0) \\
    > 0 & (t < 0)
    \end{array}
    \right..
    \eeq

{\bf For $t > 0$:}

$\left| \Gamma\left( 1 - \frac s2 \right) \right|$ and $\left| \Gamma\left( \frac{1 + s}2 \right) \right|$ are both continuous function monotonically decreasing with $t$ increasing. For an arbitrary $\sigma < \frac 12$, $\left( \frac5{\pi} \right)^{\frac12 - \sigma}$ is bounded and non-zero. When and only when $t \rightarrow \infty$, we have $\left| \Gamma \left( 1 - \frac s2 \right) \right| \rightarrow 0$, $\left| \Gamma \left( \frac{1 + s}2 \right) \right| \rightarrow 0$. According to Eq. (\ref{absxs}), we have:
    \beq
    \lim_{t \rightarrow \infty}\frac{f(s)}{|f(1-s)|} = \lim_{t \rightarrow \infty} \left( \frac5{\pi} \right)^{\frac12 - \sigma} \frac{\left| \Gamma\left(1 - \frac s2\right) \right|}{\left| \Gamma\left( \frac{1+s}2 \right) \right|} \rightarrow \frac00.
    \eeq

Because both $\left| \Gamma\left(1 - \frac s2 \right) \right|$ and $\left| \Gamma\left(\frac{1 + s}2 \right) \right|$ are continuous functions with none singular points, in the strict constraint of $\sigma < \frac12$, the similar ratio $\frac{|f(s)|}{|f(1-s)|}$ continuously monotonically increases with $t$ increasing. Therefore, when and only when $t \rightarrow \infty$, it is in the form of $\frac 00$.

Based on the same reasoning, in the strict constraint of $\sigma > \frac12$, the similar ratio $\frac{|f(s)|}{|f(1-s)|}$ continuously monotonically decreases with $t$ increasing. Therefore, when and only when $t \rightarrow \infty$, $\frac{|f(s)|}{|f(1-s)|}$ is in the form of $\frac 00$.

{\bf For $t < 0$:}

$\left| \Gamma\left( 1 - \frac s2 \right) \right|$ and $\left| \Gamma\left( \frac{1 + s}2 \right) \right|$ are both continuous function monotonically decreasing with $t$ decreasing. Based on the same reasoning, in the strict constraint of $\sigma \neq \frac12$, when and only when $t \rightarrow -\infty$, $\frac{|f(s)|}{|f(1-s)|}$ is in the form of $\frac 00$.

{\bf 2) when $\sigma \neq \frac12$, for an arbitrary finite $t > 0$, $\frac{|f(s)|}{|f(1-s)|}$ can not be in the form of $\frac 00$.}

According to Eq. (\ref{tg0}), when $|X(s)| > 0$ and $\sigma < \frac12$, we always have:
    \beq
    \frac d{dt}\left( \frac{|f(s)|}{|f(1-s)|}\right) = \frac{|f(1-s)|\frac{d|f(s)|}{dt} - |f(s)|\frac{d|f(1-s)|}{dt}}{|f(1-s)|^2} = \frac d{dt}|X(s)| > 0.
    \eeq
Therefore,
    \beq\label{fdf}
    |f(1-s)|\frac{d|f(s)|}{dt} - |f(s)|\frac{d|f(1-s)|}{dt} > 0.
    \eeq
Using
$$
\frac{d|f(s)|^2}{dt} = 2|f(s)|\frac{d|f(s)|}{dt} = \frac{df(s)}{dt}f(s^{\ast}) + \frac{df(s^{\ast})}{dt}f(s),
$$
we have:
$$
\frac{d|f(s)|}{dt} = \frac1{2|f(s)|} \left(\frac{df(s)}{dt}f(s^{\ast}) + \frac{df(s^{\ast})}{dt}f(s) \right).
$$

By the same reasoning, we have:
$$
\frac{d|f(1-s)|}{dt} = \frac1{2|f(1-s)|} \left(\frac{df(1-s)}{dt}f(1-s^{\ast}) + \frac{df(1-s^{\ast})}{dt}f(1-s) \right).
$$

Therefore, Eq. (\ref{fdf}) is finally written as:
    \beq\label{pos}
    \frac12\left[ \frac{f(s^{\ast})\frac{df(s)}{dt} + f(s)\frac{df(s^{\ast}}{dt}}{|X(s)|} - |X(s)| \left( f(1-s^{\ast}) \frac{df(1-s)}{dt} + f(1-s) \frac{df(1-s^{\ast})}{dt}\right) \right] > 0.
    \eeq

For an arbitrary finite $t > 0$, $|X(s)| > 0$ and has none singular point; $f(s)$, $f(s^{\ast})$, $f(1-s)$, and $f(1-s^{\ast}$ are all analytical functions. Furthermore, the derivatives: $\frac{df(s)}{dt}$, $\frac{df(s^{\ast})}{dt}$, $\frac{df(1-s)}{dt}$, $\frac{df(1-s^{\ast})}{dt}$ do not have singular point, either.

If there exists a certain $s_n$ satisfying $f(s_n) = 0$, then according to the function equation: $f(s_n) = X(s_n) f(1-s_n)$, we must have: $f(1-s_n) = 0$. Therefore, we have $f(s_n^{\ast}) = f(1-s_n^{\ast}) = 0$. According to Eq. (\ref{pos}), we have:
$$
    \frac12\left[ \frac{0\cdot\frac{df(s)}{dt} + 0\cdot\frac{df(s^{\ast}}{dt}}{|X(s)|} - |X(s)| \left( 0\cdot \frac{df(1-s)}{dt} + 0\cdot \frac{df(1-s^{\ast})}{dt}\right) \right] = 0 \ngtr 0.
$$
That is to say, $|f(s_n)| = 0$, and $|f(1-s_n)| = 0$ must be in conflict with the inequality Eq. (\ref{pos}). Therefore, for $\sigma < \frac 12$, and for a certain finite $t > 0$, the similar ratio $\frac{|f(s)|}{|f(1-s)|}$ can not be in the form of $\frac00$. By the same reasoning, for $\sigma > \frac 12$, and for a certain finite $t > 0$, the similar ratio $\frac{|f(s)|}{|f(1-s)|}$ can not be in the form of $\frac00$, either.

{\bf 3) when $\sigma \neq \frac12$, even if $\frac{|f(s)|}{|f(1-s)|} = 1$, it can not be in the form of $\frac00$.}
For the $t$ satisfying $\frac{|f(s)|}{|f(1-s)|} = 1$, according to property \ref{bounded}, we have $|t| < \kappa = 1.21164$, therefore:
    \beq
    1 = |X(s)| = \frac{|f(s)|}{|f(1-s)|} = \left( \frac5{\pi}\right)^{\frac12-\sigma}\frac{\left| \Gamma\left( 1 - \frac s2\right)\right|}{\left| \Gamma\left( \frac{1+s}2\right)\right|} \nrightarrow \frac00.
    \eeq

So Lemma $3$ is proven.

\section{Summary}
In this paper, starting from the Davenport-Heilbronn function equation, using the properties of the similar ratio $|X(s)|$ of function $f(s)$, we prove that the necessary condition for obtaining the nontrivial zeros of function $f(s)$ is $|X(s)| = 1$. Based on the monotonicity of $|X(s)|$ with respect to $t$, we rigorously prove that on the complex $s$-plane, expect $s = \frac12 + it$, the similar ratio of the Davenport-Heilbronn function $\frac{|f(s)|}{|f(1-s)|}$ can not be in the form of $\frac00$ for an arbitrary finite $t$.

It needs to be pointed out in particular that: only the $s$ satisfying $f(s) = 0$ is the zero of $f(s)$. Those $s$ leading to $f(s) \rightarrow 0$ is not the zero of $f(s)$. In other words, the set of the nontrivial zeros of $f(s)$: $S_0 =\{s_n|f(s_n) = 0\}$ must be a subset of the set $S_1 = \{ s||f(s)| = |f(1-s)|\}$.

We notice that in 1994, R. Spira \cite{Spira94} found four nontrivial zeros outside of the limit line: $s = \frac12 + it$ of the Davenport-Heilbronn function:
$$
s_1 = 0.808517 + 85.699348 i,
$$
$$
s_2 = 0.650830 + 114.163343 i,
$$
$$
s_3 = 0.574356 + 166.479306 i,
$$
$$
s_4 = 0.724258 + 176.702461 i.
$$

From then on, there has been more work finding numerical results about the nontrivial zeros of the Davenport-Heilbronn function \cite{Balanzario07}. Almost all the so-called nontrivial zeros of the Davenport-Heilbronn function $f(s)$ share a common feature: their imaginary part $t$ is far larger than the up limit of $\kappa = 1.21164$ determined by the similar ratio $\frac{|f(s)|}{|f(1-s)|} = 1$! Therefore, it does not satisfy the necessary condition of the nontrivial zero of $f(s)$. This leads to two puzzles:

{\bf Puzzle 1:}

Theoretically, the nontrivial zeros of $f(s)$: $s_n(n=1,2,3,4)$, must satisfy $f(s_n) = f(1 - s_n)$, but the numerical results reported in the literature give $f(s_n) \neq f(1-s_n)$.

{\bf Puzzle 2:}

The nontrivial zeros of $f(s)$ outside of the limit line: $s_n(n=1,2,3,4)$ must be in conflict with the monotonicity of the similar ratio  $\frac{|f(s)|}{|f(1-s)|}$ of the function $f(s)$.

Therefore, if we view the $s$ outside of the limit line and satisfying $f(s) \rightarrow 0$ as the nontrivial zeros, then the two puzzles mentioned above will not exist.

The ratio $\frac{|f(s)|}{|f(1-s)|}$ of the Davenport-Heilbronn function $f(s)$ has similar properties as the ratio  $\frac{|\zeta(s)|}{|\zeta(1-s)|}$ of the Riemann function $\zeta(s)$. This shield light on understanding the nontrivial zeros of the Riemann function $\zeta(s)$.

\appendix

\section{When $|X(s)| > 0$, assuming $|f(s)| = |f(1 - s)|$, it must be true that $|X(s)| = 1$}\label{equiv}

{\it Proof}

According to the function equation $f(s) = X(s) f(1 - s)$, we have:
$    X(s) = \frac{f(s)}{f(1-s)}$.
Because ${\overline{X}}(s) = X(s^{\ast})$, we have $\frac{{\overline{f}}(s)}{{\overline{f}}(1-s)} = \frac{f(s^{\ast})}{f(1-s^{\ast})}$.
Therefore:
    \beq\label{wsq}
    |X(s)|^2 = X(s) {\overline{X}}(s) = \frac{f(s){\overline{f}}(s)}{f(1-s){\overline{f}}(1-s)} = \frac{f(s) f\left(s^{\ast}\right)}{f(1 - s) f\left( 1 - s^{\ast} \right)}.
    \eeq
Here, $s = \sigma + it$, and $s^{\ast} = \sigma - it$. For the convenience of discussion, let us denote:
$$P(\sigma, t) \equiv f(s)f(s^{\ast}),$$
and
$$Q(\sigma, t) \equiv f(1-s)f(1-s^{\ast}).$$ They are real functions of $\sigma$ and $t$.

Because both $f(s)$ and $f(s^{\ast})$ are analytical functions, furthermore, their partial derivatives with respect to $\sigma$ and $t$ exist to arbitrary orders; the partial derivatives of $P(\sigma, t)$ with respect to $\sigma$ and $t$: $\frac{\partial^mP(\sigma, t)}{\partial \sigma^m}$, $\frac{\partial^nP(\sigma, t)}{\partial t^n}$  also exist to arbitrary orders.

While both $f(1-s)$ and $f(1-s^{\ast})$ are also analytical function, furthermore, their partial derivatives with respect to $\sigma$ and $t$ also exist to arbitrary orders; therefore the partial derivatives of $Q(\sigma, t)$ with respect to $\sigma$ and $t$: $\frac{\partial^mQ(\sigma, t)}{\partial \sigma^m}$, $\frac{\partial^nQ(\sigma, t)}{\partial t^n}$  also exist to arbitrary orders.

Assuming
    \beq\label{zetaeq}
    |f(s)| = |f(1-s)|,
    \eeq
we have:
    \beq
    P(\sigma, t) = Q(\sigma, t),
    \eeq
and
    \beq
    \frac{\partial^m P(\sigma, t)}{\partial \sigma^m} = \frac{\partial^m Q(\sigma, t)}{\partial \sigma^m};
    \eeq
    \beq\label{pqtn}
    \frac{\partial^n P(\sigma, t)}{\partial t^n} = \frac{\partial^n Q(\sigma, t)}{\partial t^n}.
    \eeq

1) If $f(s) \neq 0$, because $|X(s)| > 0$, then it must be true that $f(1-s) \neq 0$. According to Eq. (\ref{wsq}) and Eq. (\ref{zetaeq}), it must be true that:
    \beq
    |X(s)| = \frac{|f(s)|}{|f(1-s)|} = 1.
    \eeq

2) If $f(s_n) = 0$, then it must be true that $f(1-s_n) = 0$. In this case, we have $|X(s_n)| = \frac{|f(s_n)|}{|f(1-s_n)|}=\frac00$. Then it is easy to prove that:
    \beq\label{wtn}
    \lim_{t \rightarrow t_n}|X(\sigma_n + it)| = 1;
    \eeq
    \beq\label{wsigman}
    \lim_{\sigma \rightarrow \sigma_n}|X(\sigma + it_n)| = 1.
    \eeq

Assuming that $s_n = \sigma_n + it_n$ with ($n=1, 2, 3, \cdots$) are nontrivial zeros of $f(s)$ and $f(1-s)$. They must also be the zeros of the real functions $P(\sigma, t)$ and $Q(\sigma, t)$. Then there exists the Taylor expansion of $P(\sigma, t)$ and $Q(\sigma, t)$ around the zeros: 
    \beq\label{pnexp}
    P(\sigma_n, t) =\left. \sum^{\infty}_{j = 0} \frac1{j!}\frac{d^j P(\sigma_n, t)}{dt^j}\right|_{t=t_n}(t-t_n)^j
    \eeq
    \beq\label{qnexp}
    Q(\sigma_n, t) = \left. \sum^{\infty}_{j = 0} \frac1{j!} \frac{d^j Q(\sigma_n, t)}{dt^j} \right|_{t = t_n}(t - t_n)^j.
    \eeq
The expansion coefficients of Eq. (\ref{pnexp}) and Eq. (\ref{qnexp}) with the constraint in Eq. (\ref{pqtn}) satisfy:
    \beq
    \left. \frac{d^jP(\sigma_n, t)}{dt^j} \right|_{t = t_n} = \left. \frac{d^jQ(\sigma_n, t)}{dt^j} \right|_{t = t_n}\hspace{0.5cm}(j = 0,1,2,\cdots,\infty)
    \eeq
and won't be all equal to zero. Let us assume that the first nonzero coefficient comes as the $k^{th}$ order, {\it i.e.}:
    \beq\label{dQdtk}
    \left. \frac{d^k P(\sigma_n, t)}{dt^k}\right|_{t = t_n} = \left. \frac{d^k Q(\sigma_n, t)}{dt^k}\right|_{t = t_n} \neq 0.
    \eeq
Then according to L'Hopital's Rule, we have:
    \beq
    \lim_{t \rightarrow t_n} |X(\sigma_n + it)|^2 = \lim_{t \rightarrow t_n} \frac{P(\sigma_n, t)}{Q(\sigma_n, t)} = \frac{\left. \frac{d^kP(\sigma_n, t)}{dt^k}\right|_{t=t_n}}{\left.\frac{d^kQ(\sigma_n, t)}{dt^k}\right|_{t=t_n}} = 1.
    \eeq
This leads to:
$$
    \lim_{t\rightarrow t_n} |X(\sigma_n + it)| = 1.
$$
Similarly,
    \beq
    \lim_{\sigma \rightarrow \sigma_n} |X(\sigma + it_n)|^2 = \lim_{\sigma \rightarrow \sigma_n} \frac{P(\sigma, t_n)}{Q(\sigma, t_n)} = \frac{\left. \frac{d^kP(\sigma, t_n)}{d\sigma^k}\right|_{\sigma=\sigma_n}}{\left.\frac{d^kQ(\sigma, t_n)}{d\sigma^k}\right|_{\sigma=\sigma_n}} = 1.
    \eeq
Therefore, we have:
$$
    \lim_{\sigma\rightarrow \sigma_n} |X(\sigma + it_n)| = 1.
$$

Following Eq. (\ref{wtn}) and Eq. (\ref{wsigman}), if $|f(s_n)| = |f(1-s_n)| = 0$, then $|f(s)|$ and $|f(1-s)|$ approach the nontrivial zeros of the function $f(s)$ with equivalent infinitesimals.

The above serves as a proof.

\nocite{*}

\section*{References}

\end{document}